\documentclass{article}
\usepackage[utf8]{inputenc}
\usepackage[english]{babel}
\usepackage[utf8]{inputenc}
\usepackage{amsmath}
\usepackage{float}
\usepackage{graphicx,subfigure}
\usepackage[utf8]{inputenc}
\usepackage[english]{babel}
\usepackage{amsfonts}
\usepackage{comment} 
\usepackage{tikz}
\usepackage{mathtools}
\usepackage{pifont}
\usepackage{caption}
\usepackage[colorinlistoftodos]{todonotes}

\usepackage{amsmath,amsfonts,amssymb,amsthm,epsfig,epstopdf,array}

\newtheorem{theorem}{Theorem}
\newtheorem{lemma}[theorem]{Lemma}

\usepackage{caption}

\newcommand{\norm}[1]{\left\lVert#1\right\rVert}
\usepackage{verbatim}
\usepackage{balance}

\usepackage{enumitem}
\usepackage[bottom]{footmisc}
\usepackage[utf8]{inputenc}
\usepackage[english]{babel}

\usepackage[linkcolor=blue,colorlinks=true]{hyperref}
 \usepackage{cleveref}
 \usepackage
[
        a4paper,
        left=2cm,
        right=2cm,
        top=2cm,
        bottom=2cm,
]
{geometry}
\title{\Large \textbf{Control Barrier Functions-based Semi-Definite Programs (CBF-SDPs)}: Robust Safe Control For Dynamic Systems with Relative Degree Two Safety Indices }
\author{Jaskaran Grover$^{*}$, Changliu Liu$^{*}$,  Katia Sycara $^{*}$%
\thanks{$^{*}$J. Grover, C. Liu, K. Sycara are with the Robotics Institute,
		Carnegie Mellon University, 5000 Forbes Avenue, Pittsburgh, PA 15213, USA.
		{\tt\small \{jaskarag,cliu6,sycara\}@andrew.cmu.edu. }}%
}
\date{22 August 2022}
\begin{document}
\maketitle
\section*{Abstract}
In this draft article, we consider the problem of achieving safe control of a dynamic system for which the safety index or (control barrier function (loosely)) has relative degree equal to two. We consider parameter affine nonlinear dynamic systems and assume that the parametric uncertainty is uniform and known a-priori or being updated online through an estimator/parameter adaptation law. Under this uncertainty, the usual CBF-QP safe control approach takes the form of a robust optimization problem. Both the right hand side and left hand side of the inequality constraints depend on the unknown parameter. With the given representation of uncertainty,  the CBF-QP safe control ends up being a convex semi-infinite problem. Using two different philosophies, one based on weak duality and another based on the Lossless s-procedure, we arrive at identical SDP formulations of this robust CBF-QP problem. Thus we show that the problem of computing safe controls with known parametric uncertainty can be posed as a tractable convex problem and be solved online. (\textit{This draft is work in progress}).
\section{Introduction}
We are interested in guaranteeing safety of humans in tasks that involve space sharing interactions between a human and a robot.  We establish notation first. 
Let $\boldsymbol{x}_{A}$ denote the state of the human which could either be its position or both its position and velocity. Let $\boldsymbol{x}_{R}$ denote the state of the robot. In this work, we let $\boldsymbol{x}_{R}$ refer to the robot's position and $\boldsymbol{x}_{A}$ refer to the human's position. Further, we assume both the human and the robot have single-integrator dynamics \textit{i.e.} they are velocity controlled. The closed-loop dynamics of the human are assumed to be in a parameter-affine form as:
\begin{align}
\label{sheepdynamics}
    \Dot{\boldsymbol{x}}_{A} =\boldsymbol{G}(\boldsymbol{x}_A,\boldsymbol{x}_R)\boldsymbol{\theta} + \boldsymbol{f}(\boldsymbol{x}_A,\boldsymbol{x}_R)
\end{align}
Here  $\boldsymbol{G}(\boldsymbol{x}_A,\boldsymbol{x}_R) \in \mathbb{R}^{2 \times p}$ is a matrix of domain-specific basis functions, $\boldsymbol{\theta} \in \mathbb{R}^p$ is a parameter vector and $\boldsymbol{f}(\boldsymbol{x}_A,\boldsymbol{x}_R) \in \mathbb{R}^{2}$ is a vector of domain-specific basis functions. The dynamics of the robot are velocity controlled and can be posed as:
\begin{align}
\label{doginput}
 \Dot{\boldsymbol{x}}_{R} = \boldsymbol{u}_{R}
\end{align}

\subsection{Posing the barrier function}
Suppose there is a safety constraint on the human's motion where we require the position of the human to lie in a set described as
\begin{align}
\mathcal{S} \coloneqq \{\boldsymbol{x} \in \mathbb{R}^2 \vert h(\boldsymbol{x}) \geq 0\}
\end{align}
To ensure this condition is satisfied, we use $h(\cdot)$ as a control barrier function and evaluate its derivative over the closed-loop dynamics of human.
\begin{align}
    \dot{h} &= \nabla_{\boldsymbol{x}_A}h(\boldsymbol{x}_A)^T(\boldsymbol{G}(\boldsymbol{x}_A,\boldsymbol{x}_R)\boldsymbol{\theta} + \boldsymbol{f}(\boldsymbol{x}_A,\boldsymbol{x}_R)) \nonumber \\
    &= \boldsymbol{w}^T(\boldsymbol{x}_A,\boldsymbol{x}_R)\boldsymbol{\theta} + v(\boldsymbol{x}_A,\boldsymbol{x}_R)
\end{align}
where we have defined
\begin{align}
    \boldsymbol{w}(\boldsymbol{x}_A,\boldsymbol{x}_R) &\coloneqq \boldsymbol{G}^T(\boldsymbol{x}_A,\boldsymbol{x}_R)\nabla_{\boldsymbol{x}_A}h(\boldsymbol{x}_A) \nonumber \\
    v (\boldsymbol{x}_A,\boldsymbol{x}_R)&\coloneqq \boldsymbol{f}^T(\boldsymbol{x}_A,\boldsymbol{x}_R)\nabla_{\boldsymbol{x}_A}h(\boldsymbol{x}_A)
\end{align}
Since $\dot{h}$ does not depend on the robot control $\boldsymbol{u}_R$, we take another derivative:
\begin{align}
    \ddot{h} &= \boldsymbol{\theta}^T\bigg(\frac{\partial \boldsymbol{w}}{\partial \boldsymbol{x}_A }\dot{\boldsymbol{x}}_A + \frac{\partial \boldsymbol{w}}{\partial \boldsymbol{x}_R }\dot{\boldsymbol{x}}_R  \bigg) + \bigg(\frac{\partial v^T}{\partial \boldsymbol{x}_A }\dot{\boldsymbol{x}}_A + \frac{\partial v^T}{\partial \boldsymbol{x}_R }\dot{\boldsymbol{x}}_R  \bigg) \nonumber \\
    &= \boldsymbol{\theta}^T\bigg(\frac{\partial \boldsymbol{w}}{\partial \boldsymbol{x}_A }\dot{\boldsymbol{x}}_A + \frac{\partial \boldsymbol{w}}{\partial \boldsymbol{x}_R }\boldsymbol{u}_R  \bigg) + \bigg(\frac{\partial v^T}{\partial \boldsymbol{x}_A }\dot{\boldsymbol{x}}_A + \frac{\partial v^T}{\partial \boldsymbol{x}_R }\boldsymbol{u}_R  \bigg) \nonumber \\
    &=\bigg( \boldsymbol{\theta}^T\frac{\partial \boldsymbol{w}}{\partial \boldsymbol{x}_R } +  \frac{\partial v^T}{\partial \boldsymbol{x}_R } \bigg)\boldsymbol{u}_R + \bigg( \boldsymbol{\theta}^T\frac{\partial \boldsymbol{w}}{\partial \boldsymbol{x}_A } +  \frac{\partial v^T}{\partial \boldsymbol{x}_A}\bigg)\dot{\boldsymbol{x}}_A \nonumber \\
    &=\bigg( \boldsymbol{\theta}^T\frac{\partial \boldsymbol{w}}{\partial \boldsymbol{x}_R } +  \frac{\partial v^T}{\partial \boldsymbol{x}_R } \bigg)\boldsymbol{u}_R 
    + \bigg( \boldsymbol{\theta}^T\frac{\partial \boldsymbol{w}}{\partial \boldsymbol{x}_A } +  \frac{\partial v^T}{\partial \boldsymbol{x}_A}\bigg)
    (\boldsymbol{G}(\boldsymbol{x}_A,\boldsymbol{x}_R)\boldsymbol{\theta} 
    + \boldsymbol{f}(\boldsymbol{x}_A,\boldsymbol{x}_R)) \nonumber \\
&=\bigg( \boldsymbol{\theta}^T\frac{\partial \boldsymbol{w}}{\partial \boldsymbol{x}_R } +  \frac{\partial v^T}{\partial \boldsymbol{x}_R } \bigg)\boldsymbol{u}_R 
+ \boldsymbol{\theta}^T\frac{\partial \boldsymbol{w}}{\partial \boldsymbol{x}_A}\boldsymbol{G}(\boldsymbol{x}_A,\boldsymbol{x}_R)\boldsymbol{\theta} + 
\bigg(\frac{\partial \boldsymbol{w}}{\partial \boldsymbol{x}_A}\boldsymbol{f}(\boldsymbol{x}_A,\boldsymbol{x}_R) +\boldsymbol{G}^T(\boldsymbol{x}_A,\boldsymbol{x}_R) \frac{\partial v}{\partial \boldsymbol{x}_A} \bigg)^T \boldsymbol{\theta} + \frac{\partial v^T}{\partial \boldsymbol{x}_A} \boldsymbol{f}(\boldsymbol{x}_A,\boldsymbol{x}_R) 
\end{align}
For ensuring forward invariance of the set $\mathcal{S}$, we impose $\ddot{h} + \textcolor{red}{\alpha} \dot{h} + \textcolor{red}{\beta} h \geq 0 $ where $\textcolor{red}{\alpha}$ and $\textcolor{red}{\beta}$ are suitably chosen constants. This constraint can be reformulated as 
\begin{align}
    \ddot{h} + \textcolor{red}{\alpha} \dot{h} + \textcolor{red}{\beta} h &\geq 0 \iff \nonumber \\
     -\bigg( \boldsymbol{\theta}^T\frac{\partial \boldsymbol{w}}{\partial \boldsymbol{x}_R } +  \frac{\partial v^T}{\partial \boldsymbol{x}_R } \bigg)\boldsymbol{u}_R &\leq \boldsymbol{\theta}^T\bigg(\frac{\partial \boldsymbol{w}}{\partial \boldsymbol{x}_A}\boldsymbol{G}(\boldsymbol{x}_A,\boldsymbol{x}_R)\bigg)\boldsymbol{\theta}  
+\bigg(\frac{\partial \boldsymbol{w}}{\partial \boldsymbol{x}_A}\boldsymbol{f}(\boldsymbol{x}_A,\boldsymbol{x}_R) +\boldsymbol{G}^T(\boldsymbol{x}_A,\boldsymbol{x}_R) \frac{\partial v}{\partial \boldsymbol{x}_A} + \textcolor{red}{\alpha}\boldsymbol{w}(\boldsymbol{x}_A,\boldsymbol{x}_R) \bigg)^T \boldsymbol{\theta} + \nonumber \\
&+ \bigg(\frac{\partial v^T}{\partial \boldsymbol{x}_A} \boldsymbol{f}(\boldsymbol{x}_A,\boldsymbol{x}_R) + \textcolor{red}{\alpha} v(\boldsymbol{x}_A,\boldsymbol{x}_R)+ \textcolor{red}{\beta}h(\boldsymbol{x}_A)\bigg)
\end{align}
Define the following :
\begin{align}
\label{cdhfg}
    C &\coloneqq-\frac{\partial \boldsymbol{w}}{\partial \boldsymbol{x}_R } \nonumber \\
    \boldsymbol{d}&\coloneqq-\frac{\partial v}{\partial \boldsymbol{x}_R } \nonumber \\
    H&\coloneqq\frac{\partial \boldsymbol{w}}{\partial \boldsymbol{x}_A}\boldsymbol{G}(\boldsymbol{x}_A,\boldsymbol{x}_R) \nonumber \\
    \boldsymbol{f}&\coloneqq\frac{\partial \boldsymbol{w}}{\partial \boldsymbol{x}_A}\boldsymbol{f}(\boldsymbol{x}_A,\boldsymbol{x}_R) +\boldsymbol{G}^T(\boldsymbol{x}_A,\boldsymbol{x}_R) \frac{\partial v}{\partial \boldsymbol{x}_A} + \textcolor{red}{\alpha}\boldsymbol{w}(\boldsymbol{x}_A,\boldsymbol{x}_R) \nonumber \\
    g&\coloneqq\frac{\partial v^T}{\partial \boldsymbol{x}_A} \boldsymbol{f}(\boldsymbol{x}_A,\boldsymbol{x}_R) + \textcolor{red}{\alpha} v(\boldsymbol{x}_A,\boldsymbol{x}_R) +\textcolor{red}{\beta}h(\boldsymbol{x}_A)
\end{align}
Thus, the constraint becomes:
\begin{align}
    (\boldsymbol{\theta}^TC + \boldsymbol{d})\boldsymbol{u}_R \leq \boldsymbol{\theta}^TH\boldsymbol{\theta} + \boldsymbol{f}^T\boldsymbol{\theta} + g
\end{align}
To synthesize the safe control, we pose the following QP:
\begin{align}
\label{dog_control}
\begin{aligned}
\boldsymbol{u}^*_R&= \underset{\boldsymbol{u}_R}{\arg\min}
& & \norm{\boldsymbol{u}_R}^2 \\
& \text{subject to}
& & (\boldsymbol{\theta}^TC + \boldsymbol{d})\boldsymbol{u}_R \leq \boldsymbol{\theta}^TH\boldsymbol{\theta} + \boldsymbol{f}^T\boldsymbol{\theta} + g
\end{aligned}
\end{align}
Suppose instead of knowing the true parameter  $\boldsymbol{\theta}$, we only know a bound in which it is contained as reported by an estimator or a parameter adaptation law. For simplicity, we assume that we have an estimate $\boldsymbol{\hat{\theta}}$ and a bound $\eta$ capturing the maximum deviation of the estimate $\boldsymbol{\hat{\theta}}$ from $\boldsymbol{\theta}$. Thus, we would like to solve \eqref{dog_control} for all possible parameters that lie within $\eta$ of the estimate $\boldsymbol{\hat{\theta}}$. If using a parameter adaptation law, $\eta$ will be updated online to ensure that it $\eta(t)$ is non-increasing \textit{i.e.} the parameter uncertainty doesn't increase with time. At a given time $t=\tau$, for a fixed $\eta(\tau)$, we would like to solve the following robust QP:
\begin{align}
\label{robust QP}
\begin{aligned}
\boldsymbol{u}^*_R&= \underset{\boldsymbol{u}_R}{\arg\min}
& & \norm{\boldsymbol{u}_R}^2 \\
& \text{subject to}
& & (\boldsymbol{\tilde{\theta}}^TC + \boldsymbol{d})\boldsymbol{u}_R \leq \boldsymbol{\tilde{\theta}}^TH\boldsymbol{\tilde{\theta}} + \boldsymbol{f}^T\boldsymbol{\tilde{\theta}} + g \hspace{0.4cm} \forall \norm{\boldsymbol{\tilde{\theta}} - \boldsymbol{\hat{\theta}}}_2 \leq \eta
\end{aligned}
\end{align}
This is a robust optimization problem. We can rewrite the constraint using a new variable $\boldsymbol{z}$ by making the following substitution in \eqref{robust QP}.
\begin{align}
    \boldsymbol{\tilde{\theta}} = \boldsymbol{\hat{\theta}} + \eta \boldsymbol{z} \mbox{ where } \norm{\boldsymbol{z}}_2 \leq 1.
\end{align}
Then substituting this in the QP, we get
\begin{align}
\label{normalized QP}
\begin{aligned}
\boldsymbol{u}^*_R&= \underset{\boldsymbol{u}_R}{\arg\min}
& & \norm{\boldsymbol{u}_R}^2 \\
& \text{subject to}
& & (\boldsymbol{z}^T\tilde{C} + \boldsymbol{\tilde{d}})\boldsymbol{u}_R \leq \boldsymbol{z}^T\tilde{H}\boldsymbol{z} + \boldsymbol{\tilde{f}}^T\boldsymbol{z} + \tilde{g} \hspace{0.4cm} \forall \norm{\boldsymbol{z}}_2 \leq 1
\end{aligned}
\end{align}
Here, 
\begin{align}
    \tilde{C} &\coloneqq \eta C\nonumber \\
    \boldsymbol{\tilde{d}} &\coloneqq\boldsymbol{\hat{\theta}}^TC + \boldsymbol{d}\nonumber \\
    \tilde{H} &\coloneqq \eta^2 H\nonumber \\
    \boldsymbol{\tilde{f}} &\coloneqq \eta \boldsymbol{f} + 2\eta H \boldsymbol{\hat{\theta}}\nonumber \\
     \tilde{g} &\coloneqq \boldsymbol{\hat{\theta}}^T H \boldsymbol{\hat{\theta}}+\boldsymbol{f}^T\boldsymbol{\hat{\theta}} + g ,
\end{align}
where $C$, $\boldsymbol{d}$, $H$, $\boldsymbol{f}$ and $g$ were defined in \eqref{cdhfg}.
\subsection{Posing \ref{normalized QP} as an SDP}
\subsubsection{Approach using weak duality}
Let us recall problem \eqref{normalized QP}:
\begin{align}
\label{normalized QP recall}
\begin{aligned}
\boldsymbol{u}^*_R&= \underset{\boldsymbol{u}_R}{\arg\min}
& & \norm{\boldsymbol{u}_R}^2 \\
& \text{subject to}
& & (\boldsymbol{z}^T\tilde{C} + \boldsymbol{\tilde{d}})\boldsymbol{u}_R \leq \boldsymbol{z}^T\tilde{H}\boldsymbol{z} + \boldsymbol{\tilde{f}}^T\boldsymbol{z} + \tilde{g} \hspace{0.4cm} \forall \norm{\boldsymbol{z}}_2 \leq 1
\end{aligned}
\end{align}
The inequality constraint effectively consists of an infinite number of constraints corresponding to each $\boldsymbol{z}$ in the unit ball. These infinite constraints make the overall problem a convex semi-infinite optimization problem. A conservative approach to solve this problem is by only considering the worst constraint among those infinite constraints. This conservative problem can be posed as follows:
\begin{align}
\label{normalized QP conservative}
\begin{aligned}
\boldsymbol{u}^*_R&= \underset{\boldsymbol{u}_R}{\arg\min}
& & \norm{\boldsymbol{u}_R}^2 \\
& \text{subject to}
& &  0 \leq \inf_{\norm{\boldsymbol{z}}_2 \leq 1} \underbrace{\boldsymbol{z}^T\tilde{H}\boldsymbol{z} + (\boldsymbol{\tilde{f}}-\tilde{C}\boldsymbol{u}_R )^T\boldsymbol{z} + \tilde{g}-\boldsymbol{\tilde{d}}\boldsymbol{u}_R}_{f_0(\boldsymbol{z})}   
\end{aligned}
\end{align}
We have converted the infinite constraints into a single constraint which in itself is an optimization problem. Let's consider this inner problem now.
\begin{align}
\label{inner problem}
\begin{aligned}
&\text{inf.}
& & \boldsymbol{z}^T\tilde{H}\boldsymbol{z} + (\boldsymbol{\tilde{f}}-\tilde{C}\boldsymbol{u}_R )^T\boldsymbol{z} + \tilde{g}-\boldsymbol{\tilde{d}}\boldsymbol{u}_R  \\
& \text{subject to}
& &  \norm{\boldsymbol{z}}_2 \leq 1  
\end{aligned}
\end{align}
Note that problem  \eqref{inner problem} is equivalent to the following (possibly non-convex) QCQP:
\begin{align}
\label{inner problem QCQP}
\begin{aligned}
&\text{inf.}
& & \boldsymbol{z}^T\tilde{H}\boldsymbol{z} + (\boldsymbol{\tilde{f}}-\tilde{C}\boldsymbol{u}_R )^T\boldsymbol{z} + \tilde{g}-\boldsymbol{\tilde{d}}\boldsymbol{u}_R  \\
& \text{subject to}
& &   \boldsymbol{z}^T\boldsymbol{z} \leq 1   
\end{aligned}
\end{align}
Let $\boldsymbol{z}^*(\boldsymbol{u}_R)$ be the optimizer of \eqref{inner problem QCQP} (assume it exists). We require $\boldsymbol{z}^*(\boldsymbol{u}_R)$ to satisfy $f_0(\boldsymbol{z}^*(\boldsymbol{u}_R))\geq 0$ per the robustness requirement in \eqref{normalized QP conservative}. Thus, if we can find any lower bound for $f_0(\boldsymbol{z}^*(\boldsymbol{u}_R))$ and constrain that lower bound to be non-negative, then we will ensure $f_0(\boldsymbol{z}^*(\boldsymbol{u}_R))\geq 0$ as required by the robustness constraint. We will use weak duality to obtain a lower bound for $f_0(\boldsymbol{z}^*(\boldsymbol{u}_R))$. The Lagrangian for \eqref{inner problem QCQP} is 
\begin{align}
\label{lagrangian}
    \mathcal{L}(\boldsymbol{z},\lambda) &\coloneqq  \boldsymbol{z}^T\tilde{H}\boldsymbol{z} + (\boldsymbol{\tilde{f}}-\tilde{C}\boldsymbol{u}_R )^T\boldsymbol{z} + \tilde{g}-\boldsymbol{\tilde{d}}\boldsymbol{u}_R + \lambda(\boldsymbol{z}^T\boldsymbol{z} -1)  \nonumber \\
    & = \boldsymbol{z}^T(\tilde{H}+\lambda I)\boldsymbol{z} + (\boldsymbol{\tilde{f}}-\tilde{C}\boldsymbol{u}_R )^T\boldsymbol{z} + \tilde{g}-\boldsymbol{\tilde{d}}\boldsymbol{u}_R - \lambda 
\end{align}
Define the Lagrange dual function $g(\lambda)$ as
\begin{align}
\label{lagrange dual function}
    g(\lambda) &\coloneqq \inf_{\boldsymbol{z}\in \mathbb{R}^p} \mathcal{L}(\boldsymbol{z},\lambda)
\end{align}
If $\lambda \geq 0$, then from the lower bound property, we have that $g(\lambda) \leq f_0(\boldsymbol{z}^*(\boldsymbol{u}_R))$ \textit{i.e.} if $\lambda \geq 0$, then $g(\lambda)$ is a candidate lower bound for $f_0(\boldsymbol{z}^*(\boldsymbol{u}_R))$. Thus, the robustness constraint in \eqref{robust QP} can be posed by requiring $g(\lambda) \geq 0$. To impose this constraint, we therefore derive an expression for $g(\lambda)$ in \eqref{lagrange dual function} by using the Lagrangian defined in \eqref{lagrangian}.
\begin{align}
    g(\lambda) &\coloneqq \inf_{\boldsymbol{z}\in \mathbb{R}^p} \boldsymbol{z}^T(\tilde{H}+\lambda I)\boldsymbol{z} + (\boldsymbol{\tilde{f}}-\tilde{C}\boldsymbol{u}_R )^T\boldsymbol{z} + \tilde{g}-\boldsymbol{\tilde{d}}\boldsymbol{u}_R - \lambda \nonumber \\
    &= \begin{cases}
  -\infty  & \text{if }\tilde{H}+\lambda I \nsucceq 0   \\ 
   -\infty  &\text{if }  \tilde{H}+\lambda I \succcurlyeq 0 \mbox{ and } \boldsymbol{\tilde{f}}-\tilde{C}\boldsymbol{u}_R \notin \mathcal{R}(\tilde{H}+\lambda I) \\ 
  -\frac{1}{4} (\boldsymbol{\tilde{f}}-\tilde{C}\boldsymbol{u}_R )^T(\tilde{H}+\lambda I)^{\dagger}(\boldsymbol{\tilde{f}}-\tilde{C}\boldsymbol{u}_R ) + \tilde{g}-\boldsymbol{\tilde{d}}\boldsymbol{u}_R - \lambda &\mbox{otherwise}
\end{cases}
\end{align}
Thus, to ensure that $g(\lambda) \geq 0$ when $\lambda \geq 0$,  we require.
\begin{enumerate}
    \item $-\frac{1}{4} (\boldsymbol{\tilde{f}}-\tilde{C}\boldsymbol{u}_R )^T(\tilde{H}+\lambda I)^{\dagger}(\boldsymbol{\tilde{f}}-\tilde{C}\boldsymbol{u}_R ) + \tilde{g}-\boldsymbol{\tilde{d}}\boldsymbol{u}_R - \lambda \geq 0$
    \item  $\tilde{H}+\lambda I \succcurlyeq 0$ 
    \item  $\boldsymbol{\tilde{f}}-\tilde{C}\boldsymbol{u}_R \in \mathcal{R}(\tilde{H}+\lambda I)$
\end{enumerate}
From the Schur Complement theorem, these three conditions are equivalent to
\begin{align}
\label{lmi}
    \left[ \begin{matrix}
        \tilde{H}+\lambda I && \frac{1}{2}(\boldsymbol{\tilde{f}} - \tilde{C}\boldsymbol{u}_R)\\
        \frac{1}{2}(\boldsymbol{\tilde{f}} - \tilde{C}\boldsymbol{u}_R)^T && \tilde{g}-\boldsymbol{\tilde{d}}\boldsymbol{u}_R -\lambda
    \end{matrix}\right] \succcurlyeq 0  
\end{align}
This is a linear matrix inequality type constraint. In conjunction with $\lambda \geq 0$ and \eqref{lmi}, problem \eqref{normalized QP conservative} becomes
\begin{align}
\label{sdp weak duality}
\begin{aligned}
\boldsymbol{u}^*_R,\lambda^*= \underset{\boldsymbol{u}_R,\lambda}{\arg\min}
&  \norm{\boldsymbol{u}_R}^2 \\
 \text{subject to}
&  \left[ \begin{matrix}
        \tilde{H}+\lambda I && \frac{1}{2}(\boldsymbol{\tilde{f}} - \tilde{C}\boldsymbol{u}_R)\\
        \frac{1}{2}(\boldsymbol{\tilde{f}} - \tilde{C}\boldsymbol{u}_R)^T && \tilde{g}-\boldsymbol{\tilde{d}}\boldsymbol{u}_R -\lambda
    \end{matrix}\right] \succcurlyeq 0 \\
    &\lambda \geq 0
\end{aligned}
\end{align}
Using the epigraph trick and Schur complement theorem, we can write this as
\begin{align}
\label{singleagentsdpfinal weak duality}
\begin{aligned}
\boldsymbol{u}^*_R,\lambda^*,t^*= \underset{\boldsymbol{u}_R,\lambda,t}{\arg\min}
\hspace{0.3cm} &t \\
 \text{subject to}
&  \left[ \begin{matrix}
        \tilde{H}+\lambda I && \frac{1}{2}(\boldsymbol{\tilde{f}} - \tilde{C}\boldsymbol{u}_R)\\
        \frac{1}{2}(\boldsymbol{\tilde{f}} - \tilde{C}\boldsymbol{u}_R)^T && \tilde{g}-\boldsymbol{\tilde{d}}\boldsymbol{u}_R -\lambda
    \end{matrix}\right] \succcurlyeq 0 \\
    &\lambda \geq 0 \\
    &\left[ \begin{matrix}
        I &&  \boldsymbol{u}_R\\
       \boldsymbol{u}_R^T && t
    \end{matrix}\right] \succcurlyeq 0  
\end{aligned}
\end{align}
This is an SDP which can be easily solved online. Next, we derive the same SDP using the s-procedure.
\subsubsection{Approach using the s-procedure}
Let us recall problem \eqref{normalized QP} again:
\begin{align}
\label{normalized QP recall2}
\begin{aligned}
\boldsymbol{u}^*_R&= \underset{\boldsymbol{u}_R}{\arg\min}
& & \norm{\boldsymbol{u}_R}^2 \\
& \text{subject to}
& & (\boldsymbol{z}^T\tilde{C} + \boldsymbol{\tilde{d}})\boldsymbol{u}_R \leq \boldsymbol{z}^T\tilde{H}\boldsymbol{z} + \boldsymbol{\tilde{f}}^T\boldsymbol{z} + \tilde{g} \hspace{0.4cm} \forall \norm{\boldsymbol{z}}_2 \leq 1
\end{aligned}
\end{align}
We write the constraint $\norm{\boldsymbol{z}}_2\leq 1$ constraint as follows:
\begin{align}
    \left[ \begin{matrix}
        \boldsymbol{z} \\
        1
    \end{matrix}\right]^T \left[ \begin{matrix}
        -I & 0 \\
        0 & 1
    \end{matrix}\right] \left[ \begin{matrix}
        \boldsymbol{z} \\
        1
    \end{matrix}\right] \geq  0  \iff \boldsymbol{\tilde{z}}^TP\boldsymbol{\tilde{z}} \geq 0
\end{align}
where we have defined 
\begin{align}
    \boldsymbol{\tilde{z}} &\coloneqq (\boldsymbol{z}^T,1)^T \nonumber \\
    P &\coloneqq \left[ \begin{matrix}
        -I && 0 \\
        0 && 1
    \end{matrix}\right]
\end{align}
Then we write $(\boldsymbol{z}^T\tilde{C} + \boldsymbol{\tilde{d}})\boldsymbol{u}_R \leq \boldsymbol{z}^T\tilde{H}\boldsymbol{z} + \boldsymbol{\tilde{f}}^T\boldsymbol{z} + \tilde{g}$ as follows
\begin{align}
    \left[ \begin{matrix}
        \boldsymbol{z} \\
        1
    \end{matrix}\right]^T \left[ \begin{matrix}
        \tilde{H}& \frac{1}{2}(\boldsymbol{\tilde{f}} - \tilde{C}\boldsymbol{u}_R)\\
        \frac{1}{2}(\boldsymbol{\tilde{f}} - \tilde{C}\boldsymbol{u}_R)^T & \tilde{g}-\boldsymbol{\tilde{d}}\boldsymbol{u}_R
    \end{matrix}\right] \left[ \begin{matrix}
        \boldsymbol{z} \\
        1
    \end{matrix}\right] \geq 0 \iff \boldsymbol{\tilde{z}}^TQ\boldsymbol{\tilde{z}} \geq 0,
\end{align}
where we have defined  
\begin{align}
    Q &\coloneqq \left[ \begin{matrix}
        \tilde{H} && \frac{1}{2}(\boldsymbol{\tilde{f}} - \tilde{C}\boldsymbol{u}_R)\\
        \frac{1}{2}(\boldsymbol{\tilde{f}} - \tilde{C}\boldsymbol{u}_R)^T && \tilde{g}-\boldsymbol{\tilde{d}}\boldsymbol{u}_R
    \end{matrix}\right]
\end{align}
Thus we want $\boldsymbol{\tilde{z}}^TP\boldsymbol{\tilde{z}} \geq  0$ $\implies \boldsymbol{\tilde{z}}^TQ\boldsymbol{\tilde{z}} \geq 0$. 
Recall the statement of s-lemma \footnote{\url{https://web.stanford.edu/class/msande314/lecture15.pdf}}
\begin{lemma}
    Let $M$ and $N$ be two symmetric matrices such that there exists $\boldsymbol{u}^0$ satisfying $(\boldsymbol{u}^0)^TM\boldsymbol{u}^0 >0$. Then the implication $\boldsymbol{u}^TM\boldsymbol{u} \geq 0 \implies \boldsymbol{u}^TN\boldsymbol{u} $ holds true if and only if there exists $\lambda \geq 0$ such that $N \succcurlyeq \lambda M$.
\end{lemma}
Thus, by using $M \coloneqq P$ and $N \coloneqq Q$ in the s-lemma, $\boldsymbol{\tilde{z}}^TP\boldsymbol{\tilde{z}} \geq  0$ $\implies \boldsymbol{\tilde{z}}^TQ\boldsymbol{\tilde{z}} \geq 0$ occurs when $\exists \lambda \geq 0$ for which $Q \succcurlyeq \lambda P$ \textit{i.e.},
\begin{align}
    & \left[ \begin{matrix}
        \tilde{H} && \frac{1}{2}(\boldsymbol{\tilde{f}} - \tilde{C}\boldsymbol{u}_R)\\
        \frac{1}{2}(\boldsymbol{\tilde{f}} - \tilde{C}\boldsymbol{u}_R)^T && \tilde{g}-\boldsymbol{\tilde{d}}\boldsymbol{u}_R
    \end{matrix}\right] \succcurlyeq \lambda \left[ \begin{matrix}
        -I && 0 \\
        0 && 1
    \end{matrix}\right] 
    \iff  \left[ \begin{matrix}
        \tilde{H}+\lambda I && \frac{1}{2}(\boldsymbol{\tilde{f}} - \tilde{C}\boldsymbol{u}_R)\\
        \frac{1}{2}(\boldsymbol{\tilde{f}} - \tilde{C}\boldsymbol{u}_R)^T && \tilde{g}-\boldsymbol{\tilde{d}}\boldsymbol{u}_R -\lambda
    \end{matrix}\right] \succcurlyeq 0
\end{align}
Thus, \eqref{normalized QP} becomes 
\begin{align}
\label{singleagentsdpprefinal}
\begin{aligned}
\boldsymbol{u}^*_R,\lambda^*= \underset{\boldsymbol{u}_R,\lambda}{\arg\min}
&  \norm{\boldsymbol{u}_R}^2 \\
 \text{subject to}
&  \left[ \begin{matrix}
        \tilde{H}+\lambda I && \frac{1}{2}(\boldsymbol{\tilde{f}} - \tilde{C}\boldsymbol{u}_R)\\
        \frac{1}{2}(\boldsymbol{\tilde{f}} - \tilde{C}\boldsymbol{u}_R)^T && \tilde{g}-\boldsymbol{\tilde{d}}\boldsymbol{u}_R -\lambda
    \end{matrix}\right] \succcurlyeq 0 \\
    &\lambda \geq 0
\end{aligned}
\end{align}
Using the epigraph trick and Schur complement theorem, we can write this as
\begin{align}
\label{singleagentsdpfinal}
\begin{aligned}
\boldsymbol{u}^*_R,\lambda^*,t^*= \underset{\boldsymbol{u}_R,\lambda,t}{\arg\min}
\hspace{0.3cm} &t \\
 \text{subject to}
&  \left[ \begin{matrix}
        \tilde{H}+\lambda I && \frac{1}{2}(\boldsymbol{\tilde{f}} - \tilde{C}\boldsymbol{u}_R)\\
        \frac{1}{2}(\boldsymbol{\tilde{f}} - \tilde{C}\boldsymbol{u}_R)^T && \tilde{g}-\boldsymbol{\tilde{d}}\boldsymbol{u}_R -\lambda
    \end{matrix}\right] \succcurlyeq 0 \\
    &\lambda \geq 0 \\
    &\left[ \begin{matrix}
        I &&  \boldsymbol{u}_R\\
       \boldsymbol{u}_R^T && t
    \end{matrix}\right] \succcurlyeq 0  
\end{aligned}
\end{align}
Note that this SDP is equivalent to \eqref{singleagentsdpfinal weak duality}.
\end{document}